\documentclass[11pt]{article}
\usepackage[english]{style}
\usepackage{latexsym,amsmath,amssymb,amsfonts,graphics}

\DeclareFontFamily{U}{rsf}{}
\DeclareFontShape{U}{rsf}{m}{n}{
  <5> <6> rsfs5 <7> <8> <9> rsfs7 <10-> rsfs10}{}
\DeclareMathAlphabet{\mathscr}{U}{rsf}{m}{n}

\DeclareMathAlphabet{\mathgth}{U}{euf}{m}{n}

\DeclareFontFamily{U}{cyr}{}
\DeclareFontShape{U}{cyr}{m}{n}{
  <5> wncyr5 <6> wncyr6 <7> wncyr7 <8> wncyr8 <9> wncyr9 <10-> wncyr10}{}
\DeclareMathAlphabet{\mathcyr}{U}{cyr}{m}{n}

\input cyracc.def

\makeatletter
\def\operator@font{\sf}
\makeatother

\newcommand{\sC}{{\mathcal C}}

\newcommand{\cM}{{\mathscr M}}
\newcommand{\cN}{{\mathscr N}}

\newcommand{\dU}{{U^{\cdot}}}

\newcommand{\FMXY}{\Phi_{X\ra Y}}

\newcommand{\fmXY}{\phi_{X\ra Y}}

\newcommand{\odd}{{\mathsf{odd}}}
\newcommand{\perf}{{\mathsf{perf}}}

\newcommand{\D}{{\mathbf D}}

\newcommand{\chk}{{\scriptscriptstyle\vee}}
\newcommand{\R}{\mathbf{R}}
\newcommand{\Ld}{\mathbf{L}}
\newcommand{\lotimes}{\stackrel{\Ld}{\otimes}}

\newcommand{\free}{{\mathsf{free}}}

\DeclareMathOperator{\Pic}{Pic}

\DeclareMathOperator{\Td}{td}
\DeclareMathOperator{\ch}{ch}

\DeclareMathOperator{\Cl}{Cl}
\DeclareMathOperator{\Br}{Br}

\DeclareMathOperator{\rk}{rk}

\newcommand{\ra}{\rightarrow}

\newcommand{\lra}{\longrightarrow}

\newcommand{\scdot}{{\,\cdot\,}}

\newcommand{\C}{\mathbf{C}}

\newcommand{\Q}{\mathbf{Q}}
\newcommand{\Z}{\mathbf{Z}}

\newcommand{\iso}{\cong}

\newcommand{\pj}{\mathbf{P}}

\newcommand{\ii}{{\mathtt i}}
\renewcommand{\phi}{\varphi}

\author{%
Andrei C\u ald\u araru\thanks{Mathematics Department, University of
Wisconsin--Madison, 480 Lincoln Drive, Madison, WI 53706--1388, USA,
{\em e-mail: }{\tt andreic@math.wisc.edu}}}

\title{Non-birational Calabi-Yau threefolds that are derived equivalent}

\date{}

\begin{document}

\maketitle

\begin{abstract}
We argue that the existence of genus one fibrations with multisections
of high degree on certain Calabi-Yau threefolds implies the existence
of pairs of such varieties that are not birational, but are derived
equivalent.  It also (likely) implies the existence of non-birational
counterexamples to the Torelli problem for Calabi-Yau threefolds.
\end{abstract}

\section*{Introduction}

\paragraph
\label{thm:main}
The purpose of this paper is to prove the following theorem.
\medskip

\noindent {\bf Theorem.}  {\it Assume that there exists a smooth,
projective Calabi-Yau threefold $X$ with the following properties:
\medskip

\begin{enumerate}
\item The Picard number of $X$ is $2$,
\[ \rho = \rk\Pic(X) = 2; \]
\item $X$ admits a non-isotrivial genus one fibration $X\ra S$ of
  fiber degree $n$;
\item the group $U(\Z_n)/\Z_2$ has at least 3 elements ($\Z_2$ acts by
  negation on $\Z_n$, and $U(\Z_n)$ denotes the units in the ring $\Z_n$);
\item $X\ra S$ has at worst isolated non-reduced fibers.
\end{enumerate}
\medskip

\noindent
Then there exist smooth, projective Calabi-Yau threefolds $X'$ and
$X''$ that are not birational and have equivalent derived categories,
$\D(X') \iso \D(X'')$.  Furthermore, $X'$ and $X''$ have isomorphic
$\Z[1/2]$-Hodge structures.}

\paragraph 
The interest in a theorem of the type described above is threefold
(pun intended).  In the first place, we are interested in
understanding the classification problem of Fourier-Mukai equivalence
classes of varieties (varieties that have equivalent derived
categories).  This topic has been a central one in recent years in
algebraic geometry, having implications to the study of moduli spaces
of sheaves, the Mori program, and mirror symmetry.  In dimensions 1
and 2 the classification is complete by work of
Mukai~\cite{MukAb},~\cite{MukK3}, Bondal-Orlov~\cite{BonOrl},
Orlov~\cite{Orl}, and Bridgeland-Maciocia~\cite{BriMac}.  It is also
known that the problem is essentially only interesting for Calabi-Yau
varieties.  This follows from results of Bondal-Orlov~\cite{BonOrl}
and Toda~\cite{Tod}.  In dimension 3, we know that birational
Calabi-Yau threefolds are Fourier-Mukai equivalent by work of
Bridgeland~\cite{BriFlops}, but no examples of non-birational
Calabi-Yau threefolds with equivalent derived categories are known.

\paragraph 
\label{que:main}
{\bf Question.}  {\em Is there a pair of Fourier-Mukai equivalent
Calabi-Yau threefolds that are not birational?}

\paragraph 
Theorem~\ref{thm:main} should be regarded as strong evidence that the
answer to Question~\ref{que:main} is ``yes.''  Indeed, plenty of
examples of varieties satisfying the conditions of
Theorem~\ref{thm:main} but with $n\leq 6$ are known, and there is no
apparent reason why $n\leq 6$ would be a limitation.  The only reason
we can not find explicit examples of varieties with the required
properties is the fact that we do not know explicit equations for
elliptic curves of degree 7 or more.

\paragraph
Certain aspects of string theory also provide motivation for trying to
understand Fourier-Mukai equivalent varieties.  Kontsevich's
Homological Mirror Symmetry conjecture~\cite{Kon} predicts an
equivalence of triangulated categories between the derived category of
$X$ and the Fukaya category of $\check{X}$, when $X$ and $\check{X}$
are a mirror pair.  If $X'$ and $X''$ have the same derived category, it
is expected that they should have the same mirror as well.  Such
examples are very interesting to study.

More generally, the derived category of $X$ is conjectured to encode
the entire B-model topological quantum field theory of the open string
theory compactified on $X$~\cite{KatSha}.  Furthermore, it has been
argued~\cite{KapRoz} that this data should also be sufficient for
determining the closed string topological quantum field theory as
well.  Thus Fourier-Mukai equivalent Calabi-Yau threefolds must yield
equivalent B-model TQFT's, and thus it is useful to have examples at
hand.

\paragraph
In a third direction, we are interested in understanding how the
Torelli principle fails for Calabi-Yau threefolds.  The Torelli
problem asks how much information about a space $X$ can be recovered
from Hodge data on $X$.  For example, a typical Torelli-type theorem
asserts that if two K3 surfaces have isomorphic $\Z$-Hodge structures,
then they are themselves isomorphic.  If two spaces of the same type
(say, Calabi-Yau threefolds) have isomorphic Hodge structures but are
not themselves isomorphic, we say that they fail the Torelli
principle.

For Calabi-Yau threefolds this is known to happen.  Szendr\H oi's
examples~\cite{Sze} are pairs of Calabi-Yau threefolds that are
deformation equivalent, have isomorphic $\Z$-Hodge structures, but are
not isomorphic.  However, the varieties in these examples are pairwise
birational, so one may ask if all failures of Torelli for Calabi-Yau
threefolds are limited to birational varieties.

Theorem~\ref{thm:main} should be viewed as evidence that Torelli might
fail even for non-birational examples.  While we are unable to prove a
statement about $\Z$-Hodge structures, we argue that the varieties
$X'$ and $X''$ constructed in the theorem have isomorphic
$\Z[1/2]$-Hodge structures, and we give a convincing argument that the
isomorphism of Hodge structures might actually extend to one over
$\Z$.  Thus $X'$ and $X''$ are very close to providing an example of a
non-birational failure of Torelli.

It is worth mentioning a similar result of Szendr\H oi~\cite{Sze2}.
He investigates an example of Aspinwall and Morrison, which yields
several Calabi-Yau threefolds with isomorphic $\Z[1/5]$-Hodge
structures.  Szendr\H oi conjectures that these threefolds are not
birational~(\cite[Conjecture 0.2]{Sze2}).  The main advantage of our
approach is that, if we could construct our example, we could actually
prove that the spaces constructed are not all birational to each
other.

\paragraph
An interesting question related to the Torelli problem was posed by
Eyal Markman.  Failures of Torelli for hyperk\"ahler manifolds were
found by Namikawa~\cite{Nam} and Markman~\cite{Mar}, which were caused
by the monodromy group of the family being smaller than the group of
Hodge isometries of the integral cohomology
lattice. Theorem~\ref{thm:main} predicts the existence of
counterexamples to classical Torelli (where the group in question is
that of Hodge isometries of the integral cohomology lattice).  It
would be interesting to know whether these examples might also yield
failures of monodromy Torelli (see Griffiths~\cite{Gri} for a
discussion of the relationship between the monodromy Torelli principle
and classical Torelli).

\paragraph
The paper is organized as follows.  We collect in
Section~\ref{sec:one} required results about Fourier-Mukai transforms
and genus one fibrations.  We also include an easy example of a
variety satisfying the conditions of Theorem~\ref{thm:main} with
$n=3$.  In Section~\ref{sec:two} we prove the existence of
Fourier-Mukai partners that are not birational.  An amusing aspect of
the construction is the fact that the proof of non-birationality is
{\em not constructive}.  Explicitly, we construct three varieties that
are Fourier-Mukai partners, and we argue that at most two of them can
be birational.  We conclude with Section~\ref{sec:three} where we
prove that $X'$ and $X''$ of Theorem~\ref{thm:main} have isomorphic
$\Z[1/2]$-Hodge structures, and we discuss what is needed to obtain an
isomorphism over $\Z$.

\paragraph {\bf Acknowledgments.} This paper grew out of conversations
with Mark Gross, Paul Aspinwall, and Igor Dolgachev.  Input on the
Torelli problem was provided by Bal\'azs Szendr\H oi and Eyal Markman.

\section{Preliminaries}
\label{sec:one}

In this section we collect several results about Fourier-Mukai
transforms and genus one fibrations.

\paragraph
Let $X$ and $Y$ be smooth, projective varieties.  Given an object $E$
of $\D(X\times Y)$, define the {\em integral transform} $\FMXY^E$ with
kernel $E$ to be the functor $\D(X) \ra \D(Y)$ given by
\[ \FMXY^E(\,-\,) = \R \pi_{Y,*}(\pi_X^*(\, -\,) \lotimes E), \]
where $\pi_X, \pi_Y$ are the projections from $X\times Y$ to $X$ and
$Y$, respectively.

If an integral transform is an equivalence, it is called a {\em
Fourier-Mukai transform}.  It is known~\cite{Orl} that all
equivalences arise as integral transforms.

\paragraph
Following Mukai, we associate to any integral transform $\Phi =
\FMXY^E$ a map on cohomology
\[ \phi = \fmXY^E:H^*(X, \Q)\ra H^*(Y,\Q) \]
by the formula
\[ \fmXY^E(\,-\,) = \pi_{Y,*}(\pi_X^*(\, -\,).v(E)). \]
Here $v(E)\in H^*(X\times Y, \Q)$ is the {\em Mukai vector} of $E$,
\[ v(E) = \ch(E).\sqrt{\Td_X}. \]
This association is functorial, in the sense that the map associated
to the identity functor is the identity on $H^*(X, \Q)$, and the map
associated to $\Psi \circ \Phi$ is $\psi \circ \phi$, if $\psi$ is
associated to $\Psi$ and $\phi$ to $\Phi$.

\paragraph
\label{subsec:HHdec}
Extending this construction to complex cohomology, the map $\phi$
associated to an integral transform preserves the {\em Hochschild
grading}.  Explicitly, for every integer $i$, $\phi$ restricts to a map
\[ \phi:\bigoplus_{q-p = i} H^{p,q}(X) \ra \bigoplus_{q-p = i}
H^{p,q}(Y). \]

\paragraph
\label{subsec:MukIsom}
We have argued in~\cite{CalHH2} that Mukai's original pairing
from~\cite{MukK3} on the cohomology of K3 surfaces can be generalized
to arbitrary varieties.  For simplicity we restrict our attention to
Calabi-Yau varieties only, and then the {\em generalized Mukai
  pairing} on $H^*(X, \Q)$ is the map
\[ \langle\,-\, ,\,-\,\rangle : H^*(X, \Q) \otimes H^*(X, \Q) \ra \C
\]
defined by
\[ \langle v, w \rangle = \int_X v^\chk.w, \]
where, for a vector $v = \sum v_k$ with $v_k \in H^k(X, \Q)$,
\[ v^\chk = \sum \ii^k v_k. \]
(Here $\ii = \sqrt{1}$.)  Adjoint functors give rise to adjoint maps
on cohomology with respect to the Mukai pairing, and therefore if
$\Phi:\D(X)\ra \D(Y)$ is an equivalence, then $\phi$ is an isometry
between $H^*(X,\Q)$ and $H^*(Y, \Q)$, endowed with the Mukai pairings.

\paragraph
Let $X$ be a smooth, connected projective variety.  A flat projective
morphism $X\ra S$ to a smooth variety $S$ is called a {\em genus one
  fibration} if its generic fiber $X_\eta \ra \eta$ is a smooth curve
of genus one over the function field $\eta$ of $S$.  In general
$X_\eta$ is not expected to have any rational points, as $\eta$ is not
algebraically closed.  Equivalently, the map $X\ra S$ is not expected
to have any rational sections.  However, there is always a finite
field extension $\eta'\ra \eta$ over which $X_{\eta'}$ does have a
section.  The smallest degree of such an extension is called the {\em
  fiber degree} of the fibration $X\ra S$.  It is also the smallest
positive degree of the restriction of a relatively ample divisor on
$X/S$ to a smooth fiber of $X\ra S$.

Fix a relatively ample polarization of $X/S$, which will be implicit
from now on.  Consider the relative moduli space (in the sense of
Simpson~\cite{Sim}) of stable sheaves on the fibers of $X/S$ having
the same Hilbert polynomial as a line bundle of degree $k$ on a smooth
fiber of $X/S$.  It has a unique component containing a point
corresponding to a line bundle on a smooth fiber of $X/S$.  We'll
denote this component by $X^{(k)}$.  From the construction, it comes
with a natural quasi-projective morphism to $S$.  If $s\in S$ is a
closed point such that $X_s$ is smooth, then the fiber of $X^{(k)}$
over $s$ is isomorphic to the moduli space of line bundles of degree
$k$ on $X_s$.  Thus $X^{(k)}_s$ is isomorphic to $X_s$, although not
canonically (the choice of such an isomorphism depends on the choice
of a point on $X_s$, and such a choice can not be made globally on $X$
unless $X\ra S$ admits a section).

Now assume that $X$ has dimension at most 3.  It was argued by
Bridgeland and Maciocia~\cite{BriMac} that if $k$ is coprime to the degree of the
polarization on a smooth fiber, then every semistable torsion free
sheaf on a fiber of $X/S$ is stable, and therefore $X/S$ is
projective.  Furthermore, $X^{(k)}$ is smooth, a universal sheaf $E_0$
for the moduli problem considered exists on $X\times_S X^{(k)}$, and
extending $E_0$ by zero to $X\times X^{(k)}$ gives rise to a kernel
$E$ which induces a Fourier-Mukai transform
\[ \Phi_{X\ra X^{(k)}}^E :\D(X) \stackrel{\sim}{\lra} \D(X^{(k)}). \]
Thus $X^{(k)}\ra S$ is a genus one fibration, Fourier-Mukai equivalent
to $X$.  

\paragraph
For the sake of completeness, we include an example
from~\cite{CalEll}.  Let $X$ be a general hypersurface of bidegree
$(3,3)$ in $\pj^2\times\pj^2$.  It is a smooth Calabi-Yau threefold of
Picard number 2.  Projection from $X$ to any one of the two $\pj^2$'s
gives $X$ the structure of a genus one fibration of fiber degree 3.  A
similar example with fiber degree 5 can be found in [loc. cit.].

\section{A proof of non-birationality}
\label{sec:two}

In this section we prove, assuming given a space $X$ with the
properties required by Theorem~\ref{thm:main}, that there exists a
pair of Calabi-Yau threefolds $X'$ and $X''$ which are derived
equivalent but not birational.

\paragraph
The idea of the construction is to consider the various powers
$X^{(k)}$ of $X$, for $k$ coprime to $n$, the fiber degree of $X/S$.
They are all Calabi-Yau threefolds, admitting a genus one fibration
structure, and they are all Fourier-Mukai equivalent to each other.
In the general case, we'll argue that there are at least three
non-isomorphic such powers of $X$.  On the other hand, we will see
that in the birational equivalence class of $X=X^{(1)}$ there is at
most one other space that admits a genus one fibration.  Therefore we
conclude that at least two of the spaces $X^{(k)}$ are not birational
to each other.

There is one special case to consider, namely when all the $X^{(k)}$
spaces, admit not just one, but two genus one fibrations.  In this
case we'll argue that each one of the $X^{(k)}$ spaces is unique in
its birational equivalence class, but at least two of these spaces are
non-isomorphic.  Thus, again, we are able to find non-birational spaces
which are Fourier-Mukai equivalent.

\begin{Proposition}
\label{prop:isom}
For $k$ and $k'$ coprime to $n$, $X^{(k)}$ is isomorphic to $X^{(k')}$
as fibrations over $S$ (with induced structure from $X/S$) if and only
if $k'=\pm k \bmod n$.
\end{Proposition}

\begin{Proof}
Assume that $X^{(k)} iso X^{(k')}$ as fibrations over $S$.  The idea
of the proof is to argue that there exists an automorphism of $J/S$,
the relative Jacobian of $X/S$, that maps the Brauer class
representing $X^{(k)}$ to the Brauer class of $X^{(k')}$.  We then
argue that the only such automorphism possible is negation along the
fibers of $J/S$, which implies that $k' = \pm k \bmod n$.

Let $J/S = X^{(0)}/S$ be the relative Jacobian of $X/S$.  We will only
be interested in properties of the generic fiber of $J/S$, so we will
not be concerned with compactifying $J$.  Let $\eta$ be the generic
point of $S$, and consider $J_\eta$, the generic fiber of $J/S$.  It
is a smooth elliptic curve over the non-algebraically closed field
$\eta$.  

As argued in~\cite{CalEll}, there is a distinguished cohomology class
$\alpha$ in $\Br(J_\eta/\eta)$, the Brauer group of $J_\eta$,
corresponding to the original fibration $X/S$.  The class $\alpha$ can
be constructed either by means of Ogg-Shafarevich theory, or as an
obstruction class to the existence of a universal sheaf on
$X\times_\eta J$.  See~\cite{Cal} or~\cite{CalEll} for details.

By~\cite[6.5]{CalEll}, we can find an identification of the relative
Jacobian of $X^{(k)}/S$ with $J$, such that the Brauer class
corresponding to $X^{(k)}$ is $\alpha^k$.  This identification of the
relative Jacobian of $X^{(k)}/S$ with $J/S$ depends, however, on
knowledge of the fact that $X^{(k)}$ was constructed as the moduli
space of stable sheaves of rank 1, degree $k$ on the fibers of $X/S$.
If $X^{(k)}\iso X^{(k')}$, regarding $X^{(k)}$ as the moduli space of
line bundles of degree $k'$ on the fibers of $X/S$ may give rise to a
different identification of $J_{X^{(k)}}/S$ with $J/S$.  Thus, all we
can conclude is that $\alpha^{k'} = \phi(\alpha^{k})$ for some
automorphism $\phi$ of $J_\eta/\eta$.

The curve $J_\eta$ is an elliptic curve over the non-algebraically
closed field $\eta$.  An automorphism $\phi$ of $J_\eta/\eta$ could
either be fixing the origin of $J_\eta$, or act by translation by a
non-zero section, or a combination of both.  Consider first the case
when the origin is fixed by $\phi$.  We can regard $\phi$ as an
automorphism of the fibration $J_U/U$ over some Zariski open set $U$
in $S$.  Possibly restricting $U$ further, we can assume that all the
fibers of $J_U/U$ are smooth, and thus we get automorphisms of these
smooth elliptic curves over an algebraically closed field, fixing the
origin.  Since the fibration is not isotrivial,
by~\cite[IV.4.7]{HarAG} such an automorphism has got to be negation
along the fibers of $J_U/U$.  The corresponding action on
$\Br(J_\eta/\eta)$ maps $\alpha^k$ to $\alpha^{-k}$.  Therefore we
conclude that, with regards to automorphisms fixing the origin, the
only possibility is for $k' = \pm k \bmod n$.

We must rule out now the possibility of the existence of a translation
automorphism.  Such an automorphism must act as translation by a
non-zero section of $J_\eta/\eta$.  Such a section must be torsion:
otherwise its closure over all of $J/S$ would give a divisor which
would be linearly independent from either the zero section of $J$, or
from any divisor pulled back from $S$.  Thus we'd have $\rk \Cl J/S >
2$, contradicting the well-known fact (see, for
example,~\cite{DolGro}) that $\rk \Cl J/S = \rk \Cl X/S$, which is 2.

Assume that $J/S$ has a torsion section.  Consider a general smooth
curve $C$ in $S$, and consider the restriction $J_C\ra C$ of $J$ to
$C$.  It is a smooth elliptic fibration of dimension 2, and the
non-trivial torsion section of $J$ restricts to a non-zero torsion
section of $J_C/C$.  Its fibers are not multiple or reducible, since
the original fibration only had isolated such fibers.  It is also not
isotrivial.  But by~\cite[Proposition 5.3.4 (ii)]{CosDol} it does not
have any torsion sections, which is a contradiction.  Thus
$J/S$ does not have any torsion sections.

The reverse implication is easy and is left to the reader.
\end{Proof}

\begin{Proposition}
\label{prop:notbir}
If one of the spaces $X^{(k)}$ (which we'll call $X$ without loss of
generality) admits two structures of genus one fibration, then $X$ is
the unique smooth projective Calabi-Yau variety in its birational
equivalence class.  Otherwise, there may exist at most one more smooth
projective Calabi-Yau variety birational to $X$ admitting a genus one
fibration structure.
\end{Proposition}

\begin{Proof}
This proof is based on ideas from the Mori program.  Observe that
smooth, projective Calabi-Yau varieties are minimal in the sense of
this program.  Therefore, all other minimal varieties birational to
$X$ differ from $X$ by a finite number of flops.  

Since the rank of the Picard group of $X$ is 2, the nef cone of $X$ is
a wedge in $\R^2$ (see Figure~\ref{Fig1}).  Faces of this nef cone
correspond to extremal contractions of $X$, when they have rational
slope.  Each face could thus yield one of three cases:

\begin{itemize}
\item[--] an irrational slope face, which does not yield any geometric
  transition in the Mori program;
\item[--] a $K$-trivial fibration (a fibration with fibers which are
  either K3 or abelian surfaces, or a genus one fibration);
\item[--] a flopping contraction.
\end{itemize}

Note that there are precisely two faces of the nef cone, because of
the topology of $\R^2$.  If $X$ admits two structures of genus one
fibration, it means that both faces of the Mori cone correspond to
$K$-trivial fibrations, and therefore there are no flopping
contractions that can be performed on $X$.  Thus $X$ is the unique
minimal projective variety in its birational class.

If $X$ admits a unique structure of genus one fibration, this labels
one face of the nef cone.  The other face of the nef cone could have
irrational slope, correspond to a K3 or abelian surface fibration on
$X$, or to a flop.  In the first and second cases, we are done by an
argument similar to the previous case ($X$ is unique in its birational
class).  Otherwise, let $X'$ be the result of flopping $X$ along the
corresponding side of the nef cone.  The nef cone of $X'$ is a new
wedge in $\R^2$, adjacent to the nef cone of $X$.

\begin{figure}
\begin{center}
\resizebox{\textwidth}{!}{\includegraphics*{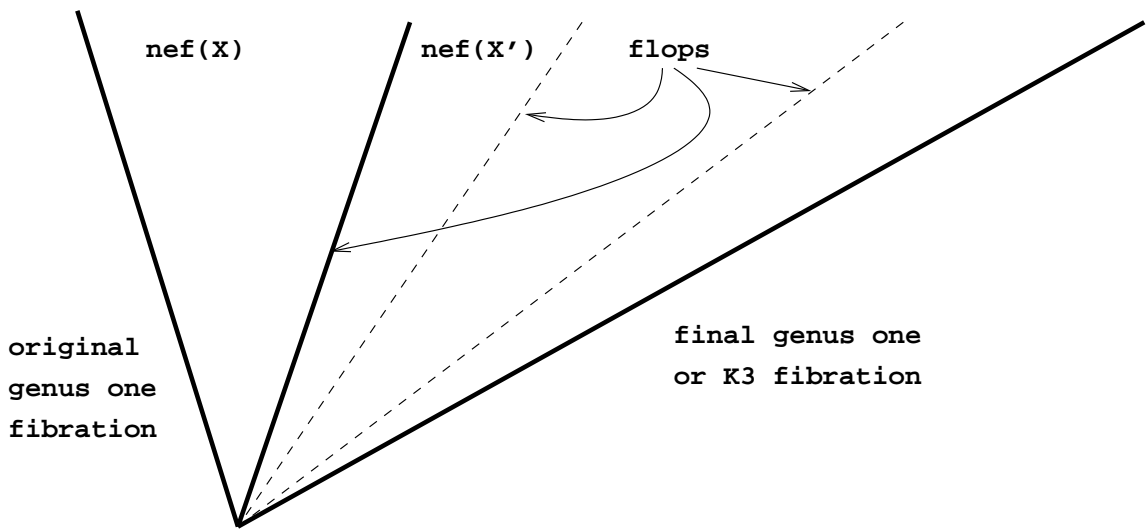}}
\caption{The nef cones of birational models of $X$}
\label{Fig1}
\end{center}
\end{figure}

Now consider the new face of the nef cone of $X'$.  It could again
have irrational slope, or yield a $K$-trivial fibration.  If it does,
our analysis of minimal models of $X$ is finished ($X$ and $X'$ are
the only minimal varieties birational to $X$, because the only
flopping contraction of $X$ yields $X'$, and vice versa).  The
$K$-trivial fiber space structure on $X'$ could be that of a genus one
fibration, or a K3 or abelian surface fibration, but in any case this
yields at most one more variety birational to $X$ with a genus one
fibration structure, and the process has stopped.

If, on the other hand, the new face of the nef cone of $X'$ yields a
flopping contraction, then $X'$ does not have any genus one fibration
structure.  Indeed, any such structure must arise from a $K$-trivial
fibration-type face of the nef cone, and both faces of
$\mathsf{nef}(X)$ yield flops, not fiber spaces.  By flopping the new
ray gives rise to a new space minimal birational model $X''$ of $X$,
but $X'$ itself does not need to be considered for examples of genus
one fibrations birational to $X$.

This process can now be continued.  At each step we consider the new
face of the nef cone of the last space constructed.  If this face
either has irrational slope, or yields a $K$-trivial fibration, we
might have arrived at one more genus one fibration on a space
birational to $X$, but the generating process of birational models of
$X$ has stopped.  Or, if this face yields a flop, we generate a new
minimal model of $X$, but the current space does not have any genus
one fibration structure.  Only the ends of the chain of nef cones can
yield minimal models of $X$ with genus one structure.  One of them is
$X$ itself, and there can be at most one other such end.
\end{Proof}

\begin{Theorem}
If a space $X$ satisfying the properties required by
Theorem~\ref{thm:main} exists, then there exist smooth, projective
Calabi-Yau threefolds $X'$ and $X''$ that are Fourier-Mukai
equivalent, but not birational.
\end{Theorem}

\begin{Proof}
Consider the collection $\sC = \{X^{(k)}\}$ for $k$ coprime to $n$,
$1\leq k \leq (n-1)/2$.   

If all the spaces in $\sC$ admit a unique structure of genus one
fibration (which must be the one arising from regarding them as
relative moduli spaces on $X/S$), then Proposition~\ref{prop:isom}
shows that they must all be distinct.  By the assumptions made on $n$,
there are at least 3 of them.  On the other hand, by
Proposition~\ref{prop:notbir}, they can not all be birational to $X$,
thus the conclusion of the theorem.

If one of the $X^{(k)}$'s in the collection admits two genus one
fibration structures, then it is unique in its birational equivalence
class.  If any other space $X^{(k')}$ in $\sC$ is non-isomorphic
to it, then $X^{(k)}$ and $X^{(k')}$ is a pair of Fourier-Mukai
equivalent, non-birational Calabi-Yau threefolds.

The only case left to consider is if all the spaces in $\sC$ are
isomorphic to each other, each admitting precisely two structures of
genus one fibration.  Note that on each space in $\sC$, there is a
distinguished such structure, arising from their construction as
moduli spaces on $X/S$.  Since there are at least 3 different
allowable values of $k$, we must be able to find $k$ and $k'$ such
that $X^{(k)}/S$ and $X^{(k')}/S$ are isomorphic as fibrations, with
the induced structure from $X/S$.  This contradicts
Proposition~\ref{prop:isom}, and therefore this case is not possible.
\end{Proof}

\section{An isomorphism of Hodge structures}
\label{sec:three}

In this section we argue first that any two Calabi-Yau threefolds that
are Fourier-Mukai equivalent have isomorphic $Z[1/2]$-Hodge
structures.  We then argue that if we can regard one of the spaces
$X'$ and $X''$ constructed in the previous section as the moduli space
of stable rank 2 vector bundles on the fibers of the other, then they
have isomorphic $\Z$-Hodge structures.

\begin{Proposition}
\label{isomCY1}
Let $X$ and $Y$ be Calabi-Yau threefolds, (complex projective manifolds, simply
connected and with trivial canonical class) such that $\D(X) \iso \D(Y)$.  Then
there exists a Hodge isometry
\[ H^3(X, \Z[1/2])_\free \iso H^3(Y, \Z[1/2])_\free, \]
where the intersection form on these groups is given by the cup product,
followed by evaluation against the fundamental class of the space, and the
subscript ``free'' denotes the torsion-free part of the corresponding group.

If $\dU$ is an element of $\D(X\times Y)$ that induces an isomorphism
$\D(X)\iso \D(Y)$, denote by $c_3'(\dU)$ the component of $c_3(\dU)$ that lies
in $H^3(X, \Z)\otimes H^3(Y, \Z)$ under the K\"unneth decomposition.  If
$c_3'(\dU)$ is divisible by 2, the isometry between the $H^3$ groups described
above is integral, i.e.\ it restricts to a Hodge isometry
\[ H^3(X, \Z)_\free \iso H^3(Y, \Z)_\free. \]
\end{Proposition}

\paragraph
It is worth observing that $H^3$ is the only cohomology group of a
Calabi-Yau threefold that carries any Hodge structure data.  Thus the
above Proposition refers to all the Hodge information of such a space.

There is, however, one more topological information that may not agree
between $X$ and $Y$, despite the above theorem.  It is the cubic form
on $H^2$, obtained by taking a class in $H^2$, cubing it, and then
integrating against the fundamental class.  

This is not too big a problem, though.  In general, when looking for
counterexamples to Torelli, we want examples that are in the same
family.  In our case, this would have to be analyzed by some other
external means, depending on the explicit geometry of the genus one
fibration we start with.  If we know that $X'$ and $X''$ are
deformation equivalent, then their cubic forms on $H^2$ will be the
same.

\paragraph
\textbf{Warning:} Throughout this section we will use a non-standard notation,
by using $H^{p,q}$ for the subgroup $H^p(X, \C) \otimes H^q(Y, \C)$ of
$H^{p+q}(X \times Y, \C)$ (from the K\"unneth decomposition), instead of the
one that comes from the Hodge decomposition.

\begin{Proof}
Using~\cite[Theorem 2.2]{Orl}, we can assume that the isomorphism
$\D(X) \iso \D(Y)$ is given by a Fourier-Mukai transform $\FMXY^\dU$
for some $\dU\in \D(X \times Y)$.  By~(\ref{subsec:MukIsom}) we get an
isometry
\[ \phi: H^*(X, \C) \iso H^*(Y, \C), \]
where the two groups are endowed with the generalized Mukai pairing.
This isometry must take $H^3(X, \C)$ to $H^3(Y,\C)$ and vice-versa.
Indeed, $H^3$ contains all the odd cohomology of a Calabi-Yau
threefold, and because the isometry respects the Hochschild grading,
it must map $H^{\odd}(X)$ to $H^{\odd}(Y)$.  It is easy to see that
the restriction of the Mukai pairing from $H^*(X, \C)$ to $H^3(X, \C)$
agrees with the usual Poincar\'e inner product (up to a constant
$-\ii$ sign).  Thus $\phi$ is an isometry $H^3(X, \C) \iso H^3(Y,
\C)$.  

Observe also that $\phi$ respects the usual Hodge decomposition of
$H^3$: each one of the terms of the Hodge decomposition is unique in
its Hochschild graded piece~(\ref{subsec:HHdec}).  Thus $\phi$ is a
Hodge isometry.

The correspondence $\phi$ is given by
\[ \phi(\scdot) = \pi_{Y, *}(\pi_X^*(\scdot).\ch(\dU).
\sqrt{\Td(X\times Y)}). \]
Write $u_i$ for $c_i(\dU)$, and $r$ for $\rk(\dU)$.  

Observe that the only part that contributes to the $H^3(Y, \C)$ component of
$\phi(x)$ (for some $x\in H^3(X, \C)$) is the $H^{6,3}$ part of
\[ \pi_X^*(x).\ch(\dU).\sqrt{\Td(X\times Y)}, \]
and thus, since $\pi_X^*(x)\in H^{3,0}$, we are only interested in the
$H^{3,3}$ component of
\[ \ch(\dU).\sqrt{\Td(X\times Y)}. \]
Now 
\[ \sqrt{\Td(X\times Y)} = 1 + \frac{1}{24} (\pi_X^* c_2(X) +\pi_Y^* c_2(Y)) + 
\mbox{ higher order terms,} 
\]
so that the $H^{3,3}$ component of 
\[ \ch(\dU).\sqrt{\Td(X\times Y)} \]
is just $\ch_3(\dU)$.  Indeed, $\pi_X^*c_2(X) \in H^{4,0}$, so it can not give
an element of $H^{3,3}$ by multiplication with anything.
We have 
\[ \ch_3(\dU) = \frac{1}{6}(c_1^3(\dU) - 3 c_1(\dU)c_2(\dU) + 3 c_3(\dU)); \]
$c_1(\dU)\in H^2(X\times Y, \C)$, and since $H^1(X, \C) = H^1(Y, \C) = 0$, it
must belong to $H^{2,0}\oplus H^{0,2}$.  Therefore $c_1^3(\dU)$ can not have
any $H^{3,3}$ component.  Similarly for $c_1(\dU)c_2(\dU)$.  Hence we conclude
that the only contribution to the map $H^3(X, \C)$ $\ra H^3(Y, \C)$ comes from
$c_3(\dU)$, in the form
\[ \phi|_{H^3(X, \C)}(\scdot) = \frac{1}{2}\pi_{Y,*}(\pi_X^*(\scdot).c_3'(\dU)
),\]
and hence the result.
\end{Proof}

\begin{Proposition}
\label{zhodge}
Just as before, let $X$ and $Y$ be Calabi-Yau threefolds with
equivalent derived categories.  Assume that there exists a codimension
2 subvariety $i:Z\hookrightarrow X\times Y$ which is integral and
locally a complete intersection, and an element $\dU\in \D_\perf(Z)$
such that if one takes $i_* \dU\in \D(X\times Y)$, then
$\FMXY^{i_*\dU}:\D(X) \ra \D(Y)$ is an equivalence of categories.
Assume furthermore that $\dU$ has rank 2 on $Z$.  Then the isometry of
Proposition~\ref{isomCY1} is integral, i.e.\ it restricts to a Hodge
isometry
\[ H^3(X, \Z)_\free \iso H^3(Y, \Z)_\free. \]
\end{Proposition}

\begin{Proof}
As seen before, to prove the integrality of the isomorphism we would need to
show that $c_3(i_*\dU)$ is divisible by 2.  We use Grothendieck-Riemann-Roch to
compute $c_3(i_*\dU)$ (we can apply it because $i$ is a locally complete
intersection projective morphism).  We have
\[ \ch(i_*\dU) = i_*(\ch(\dU).\Td(-\cN)), \]
where $\cN$ is the normal sheaf of $Z$ in $X\times Y$, and the negation is
taken in the Grothendieck group (see~\cite[Expos\'e 0]{SGA6}).

The left hand side of the above equality is
\[ r+c_1+\frac{1}{2}(c_1^2-c_2) + \frac{1}{6}(c_1^3 - 3 c_1 c_2 + 3c_3) +
\ldots, \] 
where $c_i = c_i(i_*\dU)$.  Obviously, $r=c_1=0$ because $Z$ has complex
codimension 2 (also from the right hand side of the equality), so the real
dimension 6 part of the left hand side consists of just
\[ \frac{1}{2} c_3(i_*\dU). \]

On the other hand, the real dimension 6 part of the right hand side comes from
the dimension 2 part of
\[ \ch(\dU).\Td(-\cN), \]
which consists of 
\[ c_1(\dU) - \frac{1}{2} \rk(\dU) c_1(\cN). \]
We conclude that 
\[ c_3(i_*\dU) = 2 i_*c_1(\dU) - \rk(\dU) i_* c_1(\cN). \]
The first term is obviously divisible by 2, and the second one is
divisible by 2 because we have assumed $\rk(\dU) = 2$.  The result follows.
\end{Proof}

\paragraph
We conjecture that the arguments of the previous section can be
repeated to get two spaces $X'$ and $X''$ which are Fourier-Mukai
equivalent, non-biratonal, and having the following property: one can
be regarded as a component of the moduli space of stable torsion free
sheaves of rank 2, degree $k$ on the other one, for some $k$.

If this is the case, then probably Proposition~\ref{zhodge} can be
applied to these spaces (the universal sheaf $U$ would be supported
on $X'\times_S X''$ which is codimension 2 in $X'\times X''$, and we
in fact expect $U$ to be a rank 2 vector bundle on $X'\times_S
X''$).  Then $X'$ and $X''$ would have isomorphic $\Z$-Hodge
structures, thus giving a complete example of two non-birational
Calabi-Yau varieties with isomorphic Hodge structures.

\paragraph
The reason we expect to be able to regard $X''$ as a moduli space of
rank 2 (instead of rank 1) vector bundles on the fibers of $X'$ is the
following.  Given an elliptic curve $E$, and coprime integers $r$ and
$d$ with $r>0$, then the moduli space of stable vector bundles on $E$
of rank $r$, degree $d$ is again isomorphic to $E$ (but the choice of
isomorphism depends on the choice of origin of $E$).  Thus starting
with $X$, and considering relative moduli spaces of stable sheaves of
rank 2, degree $k$ on the fibers of $X/S$ for various $k$, we get
genus one fibrations.  The relative Jacobian of these fibrations is
isomorphic to $J/S$, as these fibrations have the same fibers as
$X/S$.  We believe it is highly unlikely that these new fibrations
would be represented in $\Br(J/S)$ by anything other than a power of
$\alpha$ -- there does not appear to be any good way to create from
nothing new classes in the Brauer group.

If the above reasoning is correct, then the spaces we are considering
are again the $X^{(k)}$'s we have considered in Section~\ref{sec:two},
but now regarded as moduli spaces of rank 2 vector bundles.  This does
not yet say that any one of them can be regarded as a moduli space of
rank 2 vector bundles on any other one of them.  But we expect that
there exists a simple formula of the form
\[ X^{(f(2,k, k' ,n))} = \cM_{X^{(k)}}(0, 2, k') \]
where $f$ is a simple polynomial in $k$, $k'$ and $n$, and
$\cM_{X^{(k)}}(0,2,k')$ represents the moduli space of rank 2, degree
$k'$ stable sheaves on the fibers of $X^{(k)}$.  Then, for an
appropriate choice of $n$, one would be able to argue that {\em every}
$X^{(k)}$ is a moduli space of stable sheaves of rank 2, degree $l$ on
the fibers of {\em any} other such space.

\end{document}